\documentclass[11pt,twoside,reqno]{amsart}

\usepackage{microtype}
\usepackage[OT1]{fontenc}
\usepackage{type1cm}
\usepackage{amssymb}
\usepackage[left=2.3cm,top=3cm,right=2.3cm]{geometry}        
\usepackage[11pt]{moresize}

\numberwithin{equation}{section}
\numberwithin{figure}{section}

\theoremstyle{plain}

\newtheorem{theorem}{Theorem}[section]

\theoremstyle{remark}

\theoremstyle{definition}
\newtheorem{definition}[theorem]{Definition}

\newtheorem*{notations*}{Notations}
\newtheorem*{notation*}{Notation}

\renewcommand{\AA}{\mathcal{A}}
\newcommand{\HH}{\mathcal{H}}
\newcommand{\LL}{\mathcal{L}}
\newcommand{\PP}{\mathcal{P}}
\newcommand{\MM}{\mathcal{M}}

\newcommand{\BB}{\mathcal{B}}

\newcommand{\DD}{\mathcal{D}}
\newcommand{\R}{\mathbb{R}}

\newcommand{\N}{\mathbb{N}}

\newcommand{\TD}{\mathcal{TD}}

\newcommand{\eps}{\varepsilon}

\newcommand{\dist}{\operatorname{dist}}

\newcommand{\Tan}{\operatorname{Tan}}

\renewcommand{\epsilon}{\varepsilon}
\renewcommand{\rho}{\varrho}
\renewcommand{\phi}{\varphi}

\renewcommand{\hat}{\widehat}

\renewcommand{\iint}{\int\hspace{-0.09in}\int}

  \newcommand{\la}{\langle}
  \newcommand{\ra}{\rangle}

\DeclareMathOperator{\spt}{spt}

\DeclareMathOperator{\ldimloc}{\underline{dim}_{loc}}

\DeclareMathOperator{\dimh}{dim_H}

\DeclareMathOperator{\ldimh}{\underline{dim}_H}

\DeclareMathOperator*{\essinf}{ess\,inf}

\newcommand{\dd}{\,\mathrm{d}}

\begin{document}

\title{Scenery flow, conical densities, and rectifiability}

\author{Antti K\"aenm\"aki}
\address{Department of Mathematics and Statistics \\
         P.O.\ Box 35 (MaD) \\
         FI-40014 University of Jyv\"askyl\"a \\
         Finland}
\email{antti.kaenmaki@jyu.fi}

\let\thefootnote\relax\footnote{\emph{Date}: December 1, 2014.}

\begin{abstract}
We present an application of the recently developed ergodic theoretic machinery on scenery flows to a classical geometric measure theoretic problem in Euclidean spaces. We also review the enhancements to the theory required in our work. Our main result is a sharp version of the conical density theorem, which we reduce to a question on rectifiability.
\end{abstract}

\maketitle


\section{Introduction}

We survey a recent advance in the study of scenery flows and show how it can be applied in a classical question in geometric measure theory which a priori does not involve any dynamics. The reader is prompted to recall the expository article of Fisher \cite{Fisher2004} where it was discussed how the scenery flow is linked to rescaling on several well-studied structures, such as geodesic flows, Brownian motion, and Julia sets. The purpose of this note is to continue that line of introduction.

The idea behind the scenery flow has been examined in many occasions. Authors have considered the scenery flow for specific sets and measures arising from dynamics; see e.g.\ \cite{BedfordFisher1996, BedfordFisher1997, BedfordFisherUrbanski2002, FergusonFraserSahlsten2013, Zahle88}. Abstract scenery flows have also been studied with a view on applications to special sets and measures, again arising from dynamics or arithmetic; see e.g.\ \cite{Gavish2011, Hochman2010, HochmanShmerkin2013}. The main innovation of the recent article by K\"aenm\"aki, Sahlsten, and Shmerkin \cite{KaenmakiSahlstenShmerkin2014b} is to employ the general theory initiated by Furstenberg \cite{Furstenberg2008}, greatly developed by Hochman \cite{Hochman2010} and extended by K\"aenm\"aki, Sahlsten, and Shmerkin \cite{KaenmakiSahlstenShmerkin2014a}, to classical problems in geometric measure theory.

One of the most fundamental concepts of geometric measure theory is that of rectifiability. It is a measure-theoretical notion for smoothness and to a great extend, geometric measure theory is about studying rectifiable sets. The foundations of geometric measure theory were laid by Besicovitch \cite{Besicovitch1928, Besicovitch1929}. For various characterizations and properties of rectifiability the reader is referred to the book of Mattila \cite{Mattila1995}. In conical density results, the idea is to examine how a measure is distributed in small balls. Finding conditions that guarantee the measure to be effectively spread out in different directions is a classical question going back to Besicovitch \cite{Besicovitch1938} and Marstrand \cite{Marstrand1954}. For an account of the development on conical density results the reader is referred to the survey of K\"aenm\"aki \cite{Kaenmaki2010}.

The scenery flow is a well-suited tool to address problems concerning conical densities. The cones in question do not change under magnification and this allows to pass information between the original measure and its tangential structure. In fact, we will see that there is an intimate connection between rectifiability and conical densities.

This exposition comes in two parts. In the first part, we review dynamical aspects of the scenery flow and in the second part, we focus on geometric measure theory.

\section{Dynamics of the scenery flow}

Let $(X,\BB,P)$ be a probability space. We shall assume that $X$ is a metric space and $\BB$ is the Borel $\sigma$-algebra on $X$. Write $\R_+ = [0,\infty)$. A (one-sided) \emph{flow} is a family $(F_t)_{t \in \R_+}$ of measurable maps $F_t \colon X \to X$ for which
$$
  F_{t+t'} = F_{t} \circ F_{t'}, \quad t,t' \in \R_+.
$$
In other words, $(F_t)_{t \in \R_+}$ is an additive $\R_+$ action on $X$. We also assume that $(x,t) \mapsto F_t(x)$ is measurable.

We say that a set $A \in \BB$ is $F_t$ \emph{invariant} if $P(F_t^{-1} A \triangle A)=0$ for all $t \ge 0$. If $F_t P = P$ for all $t \ge 0$, then we say that $P$ is $F_t$ \emph{invariant}. In this case, we call $(X,\BB,P,(F_t)_{t \in \R_+})$ a \emph{measure preserving flow}. Furthermore, a measure preserving flow is \emph{ergodic}, if for all $t \geq 0$ the measure $P$ is ergodic with respect to the transformation $F_t \colon X \to X$, that is, for all $F_t$ invariant sets $A \in \BB$ we have $P(A) \in \{ 0,1 \}$.

\begin{theorem}[Birkhoff ergodic theorem]
  If $(X,\BB,P,(F_t)_{t \in \R_+})$ is an ergodic measure preserving flow, then for a $P$ integrable function $f \colon X \to \R$ we have
  \begin{equation*}
    \lim_{T \to \infty}  \frac{1}{T} \int_0^T f(F_t x) \dd t = \int f \dd P
  \end{equation*}
  for $P$-almost all $x \in X$.
\end{theorem}

We write $\omega \sim P$ to indicate that $\omega$ is chosen randomly according to the measure $P$. 

\begin{theorem}[Ergodic decomposition]
  Any $F_t$ invariant measure $P$ can be decomposed into ergodic components $P_\omega$, $\omega \sim P$, such that
  \begin{equation*}
    P = \int P_\omega \dd P(\omega).
  \end{equation*}
  This decomposition is unique up to $P$ measure zero sets.
\end{theorem}

Let us next define the scenery flow. We equip $\R^d$ with the usual Euclidean norm and the induced metric. Denote the closed unit ball by $B_1$. Let $\MM_1 := \PP(B_1)$ be the collection of all Borel probability measures on $B_1$ and $\MM_1^* := \{ \mu \in \MM_1 : 0 \in \spt(\mu) \}$. Here $\spt(\mu)$ is the support of $\mu$. To avoid any confusion, measures on measures will be called \emph{distributions}. We define the \emph{magnification} $S_t\mu$ of $\mu \in \MM_1^*$ at $0$ by setting
\begin{equation*}
  S_t\mu(A) := \frac{\mu(e^{-t}A)}{\mu(B(0,e^{-t}))}, \quad A \subset B_1.
\end{equation*}
In other words, the measure $S_t\mu$ is obtained by scaling $\mu|_{B(0,e^{-t})}$ into the unit ball and normalizing. Due to the exponential scaling, $(S_t)_{t \in \R_+}$ is a flow in the space $\MM_1^*$ and we call it the \emph{scenery flow} at $0$. An $S_t$ invariant distribution $P$ on $\MM_1^*$ is called \emph{scale invariant}. Although the action $S_t$ is discontinuous (at measures $\mu$ with $\mu(\partial B(0,r)) > 0$ for some $0<r<1$) and the set $\MM_1^* \subset \MM_1$ is not closed, we shall witness that the scenery flow behaves in a very similar way to a continuous flow on a compact metric space.

With the scenery flow we are now able to define tangent measures and distributions. Let $\mu$ be a Radon measure and $x \in \spt(\mu)$. We want to consider the scaling dynamics when magnifying around $x$. Let $T_x\mu(A) := \mu(A+x)$ and define $\mu_{x,t} := S_t(T_x\mu)$. Then the one-parameter family $(\mu_{x,t})_{t \in \R_+}$ is called the \emph{scenery flow} at $x$. Accumulation points of this scenery in $\MM_1$ will be called \emph{tangent measures} of $\mu$ at $x$ and the family of tangent measures of $\mu$ at $x$ is denoted by $\Tan(\mu,x) \subset \MM_1$. However, we are not interested in a single tangent measure, but the whole statistics of the scenery $\mu_{x,t}$ as $t \to \infty$. We remark that we have slightly deviated from Preiss' original definition of tangent measures, which corresponds to taking weak limits of unrestricted blow-ups; see \cite{Preiss1987}.

\begin{definition}[Tangent distributions]
  A \emph{tangent distribution} of $\mu$ at $x \in \spt(\mu)$ is any weak limit of
  \begin{equation*}
    \la \mu \ra_{x,T} := \frac{1}{T} \int_0^T \delta_{\mu_{x,t}} \dd t
  \end{equation*}
  as $T \to \infty$. The family of tangent distributions of $\mu$ at $x$ is denoted by $\TD(\mu,x) \subset \PP(\MM_1^*)$.
\end{definition}

If the limit above is unique, then, intuitively, it means that the collection of views $\mu_{x,t}$ will have well-defined statistics when zooming into smaller and smaller neighbourhoods of $x$. The integration above makes sense since we are on a convex subset of a topological linear space. We emphasize that tangent distributions are measures on measures. Notice that the set $\TD(\mu,x)$ is non-empty and compact at $x \in \spt(\mu)$. Moreover, the support of each $P \in \TD(\mu,x)$ is contained in $\Tan(\mu,x)$.

According to Preiss' well-known principle, tangent measures to tangent measures are tangent measures; see \cite[Theorem 2.12]{Preiss1987}. We shall define an analogous condition for distributions. We say that a distribution $P$ on $\MM_1$ is \emph{quasi-Palm} if for any Borel set $\AA \subset \MM_1$ with $P(\AA) = 1$ it holds that for $P$-almost every $\nu \in \AA$ and for $\nu$-almost every $z \in \R^d$ there exists $t_z > 0$ such that for $t \ge t_z$ we have $B(z, e^{-t}) \subset B_1$ and
\begin{equation*}
  \nu_{z,t} \in \AA.
\end{equation*}
This version of the quasi-Palm property actually requires that the unit sphere of the norm is a $C^1$ manifold and does not contain line segments; see \cite[Lemma 3.23]{KaenmakiSahlstenShmerkin2014b}. The Euclidean norm we use of course satisfies this requirement. If we were considering unrestricted blow-ups, then the requirement for $B(z,e^{-t})$ to be contained in $B_1$ could be dropped. Roughly speaking, the quasi-Palm property guarantees that the null sets of the distributions are invariant under translations to a typical point of the measure.

\begin{definition}[Fractal distributions]
  A distribution $P$ on $\MM_1$ is a \emph{fractal distribution} if it is scale invariant and quasi-Palm. A fractal distribution is an \emph{ergodic fractal distribution} if it is ergodic with respect to $S_t$.
\end{definition}

It follows from the Besicovitch density point theorem that ergodic components of a fractal distribution are ergodic fractal distributions; see \cite[Theorem 1.3]{Hochman2010}.

A general principle is that tangent objects enjoy some kind of spatial invariance. For tangent distributions, a very powerful formulation of this principle is the following theorem of Hochman \cite[Theorem 1.7]{Hochman2010}. The result is analogous to a similar phenomenon discovered by M\"orters and Preiss \cite[Theorem 1]{MortersPreiss1998}.

\begin{theorem} \label{thm:hochman}
  For any Radon measure $\mu$ and $\mu$-almost every $x$, all tangent distributions of $\mu$ at $x$ are fractal distributions.
\end{theorem}

Notice that as the action $S_t$ is discontinuous, even the scale invariance of tangent distributions or the fact that they are supported on $\MM_1^*$ are not immediate, though they are perhaps expected. The most interesting part in the above theorem is that a typical tangent distribution satisfies the quasi-Palm property.

Hochman's result is proved by using CP processes which are Markov processes on the dyadic scaling sceneries of a measure introduced by Furstenberg \cite{Furstenberg1970, Furstenberg2008}. Let $\DD$ be a partition of $[-1,1]^d$ into $2^d$ cubes of side length $1$. Given $x \in [-1,1]^d$, let $D(x)$ be the only element of $\DD$ containing it. If $D \in \DD$, then we write $T_D$ for the orientation preserving homothety mapping from $\overline{D}$ onto $[-1,1]^d$. Define the \emph{CP magnification} $M$ on $\Omega := \PP([-1,1]^d) \times [-1,1]^d$ by setting
\begin{equation*}
  M(\mu,x) := \bigl( T_{D(x)}\mu/\mu(D(x)), T_{D(x)}(x) \bigr).
\end{equation*}
This is well-defined whenever $\mu(D(x))>0$. Note that, since zooming in is done dyadically, it is important to keep track of the orbit of the point that is being zoomed upon. A distribution $Q$ on $\Omega$ is \emph{adapted} if there is a disintegration
\begin{equation*}
  \int f(\nu,x) \dd Q(\nu,x) = \iint f(\nu,x) \dd\nu(x) \dd \overline{Q}(\nu)
\end{equation*}
for all $f \in C(\Omega)$. Here $\overline{Q}$ is the projection of $Q$ onto the measure component. In other words, $Q$ is adapted if choosing a pair $(\mu,x)$ according to $Q$ can be done in a two-step process, by first choosing $\mu$ according to $\overline{Q}$ and then choosing $x$ according to $\mu$. A distribution on $\Omega$ is a \emph{CP distribution} if it is $M$ invariant and adapted.

The \emph{micromeasure distribution} of $\mu$ at $x \in \spt(\mu)$ is any weak limit of
\begin{equation*}
  \la \mu,x \ra_N := \frac{1}{N} \sum_{k=0}^{N-1} \delta_{M^k(\mu,x)}.
\end{equation*}
By compactness of $\PP(\Omega)$, the family of micromeasure distributions is non-empty and compact, and by \cite[Proposition 5.4]{Hochman2010}, each micromeasure distribution is adapted. Furthermore, if the \emph{intensity measure} of a micromeasure distribution $Q$ defined by
\begin{equation*}
  [Q](A) := \int \mu(A) \dd\overline{Q}(\mu), \quad A \subset [-1,1]^d,
\end{equation*}
is the normalized Lebesgue measure, then $Q$ is M invariant. By adaptedness, this is the case for any weak limit of $\la \mu+z,x+z \ra_N$ for Lebesgue almost all $z \in [-1/2,1/2]^d$; see \cite[Proposition 5.5(2)]{Hochman2010}. In other words, by slightly adjusting the dyadic grid, a micromeasure distribution can be seen to be a CP distribution. The family of CP distributions having Lebesgue intensity is compact; see \cite[Lemma 3.4]{KaenmakiSahlstenShmerkin2014a}.

If $Q$ is a CP distribution, then the system $(\Omega,M,Q)$ is a stationary one-sided process $(\xi_n)_{n \in \N}$ with $\xi_1 \sim Q$ and $M\xi_n = \xi_{n+1}$. Considering its two-sided extension, we see that there exists a natural extension $\hat Q$ supported on the Cartesian product of all Radon measures and $[-1,1]^d$. A \emph{centering} of $\hat Q$ is a push-down of the suspension flow of $\hat Q$ under the unrestricted magnification of $\mu$ at $x$. For a precise definition, see \cite[Definition 1.13]{Hochman2010}. By \cite[Theorem 1.14]{Hochman2010}, a centering of $\hat Q$ is an unrestricted fractal distribution. We remark that \cite{Hochman2010} and \cite{KaenmakiSahlstenShmerkin2014a} use $L^\infty$ norm to allow an easier link between CP processes and fractal distributions. By \cite[Appendix A]{KaenmakiSahlstenShmerkin2014a}, the results are independent of the choice of the norm and hence, our use of the Euclidean norm is justified.

Relying on the above, we are now able to give an outline for the proof of Theorem \ref{thm:hochman}. If $P = \lim_{k \to \infty} \la \mu \ra_{x,N_k}$ is a tangent distribution, then, passing to a subsequence, define a micromeasure distribution $Q = \lim_{i \to \infty} \la \mu,x \ra_{N_{k(i)}}$. Slightly adjusting the dyadic grid, we see that $Q$ is a CP distribution with Lebesgue intensity. Thus, by \cite[Proposition 5.5(3)]{Hochman2010}, $P$ is the restriction of the centering of $\hat Q$ and hence, $P$ is a fractal distribution.

Although fractal distributions are defined in terms of seemingly strong geometric properties, the family of fractal distributions is in fact very robust. The following theorem is due to K\"aenm\"aki, Sahlsten, and Shmerkin \cite[Theorem A]{KaenmakiSahlstenShmerkin2014a}.

\begin{theorem} \label{thm:compact}
  The family of fractal distributions is compact.
\end{theorem}

The result may appear rather surprising since the scenery flow is not continuous, its support is not closed, and, more significantly, the quasi-Palm property is not a closed property. The proof of this result is also based on the interplay between fractal distributions and CP processes. We have already seen that each CP distribution defines a fractal distribution. The converse is also true. Let us first assume that $P$ is an ergodic fractal distribution. If $f$ is a continuous function defined on $\PP(\MM_1)$, then, by the Birkhoff ergodic theorem, we have
\begin{equation*}
  \lim_{T \to \infty} \frac{1}{T} \int_0^T f(S_t\mu) \dd t = \int f \dd P
\end{equation*}
for $P$-almost all $\mu$. Considering a countable dense set of continuous functions $f$ and applying the quasi-Palm property, it follows that
\begin{equation} \label{eq:erg_usm}
  \lim_{T \to \infty} \la \mu \ra_{x,T} = P
\end{equation}
for $P$-almost all $\mu$ and for $\mu$-almost all $x$; see \cite[Theorem 3.9]{Hochman2010}. As we already have seen, any tangent distribution can be expressed as the restriction of the centering of an extended CP distribution having Lebesgue intensity. Thus, by \eqref{eq:erg_usm}, the same holds for ergodic fractal distributions. Relying on the ergodic decomposition, this observation can be extended to non-ergodic fractal distributions; see \cite[Theorem 1.15]{Hochman2010}. Therefore, since the family of CP distributions with Lebesgue intensity is compact, to prove Theorem \ref{thm:compact}, it suffices to show that the centering is a continuous operation. This is done in \cite[Lemmas 3.5 and 3.6]{KaenmakiSahlstenShmerkin2014a}.

Together with convexity and the uniqueness of the ergodic decomposition, Theorem \ref{thm:compact} implies that the family of fractal distributions is a Choquet simplex. Recall that a Poulsen simplex is a Choquet simplex in which extremal points are dense. Note that the set of extremal points is precisely the collection of ergodic fractal distributions. The following theorem is proved by K\"aenm\"aki, Sahlsten, and Shmerkin \cite[Theorem B]{KaenmakiSahlstenShmerkin2014a}.

\begin{theorem} \label{thm:poulsen}
  The family of fractal distributions is a Poulsen simplex.
\end{theorem}

The proof is again based on the interplay between fractal distributions and CP processes. We prove that ergodic CP processes are dense by constructing a dense set of distributions of random self-similar measures on the dyadic grid. This is done by first approximating a given CP process by a finite convex combination of ergodic CP processes, and then, by splicing together those finite ergodic CP processes, constructing a sequence of ergodic CP processes converging to the convex combination. Roughly speaking, splicing of measures consists in pasting together a sequence of measures along dyadic scales. Splicing is often employed to construct measures with a given property based on properties of the component measures. For details, the reader is referred to \cite[\S 4]{KaenmakiSahlstenShmerkin2014a}.

In geometric considerations, we usually construct a fractal distribution satisfying certain property. We often want to transfer that property back to a measure. This leads us to the concept of generated distributions.

\begin{definition}[Uniformly scaling measures]
  We say that a measure $\mu$ \emph{generates} a distribution $P$ at $x$ if
  \begin{equation*}
    \TD(\mu,x) = \{ P \}.
  \end{equation*}
  If $\mu$ generates $P$ for $\mu$-almost all $x$, then we say that $\mu$ is a \emph{uniformly scaling measure}.
\end{definition}

One can think that the uniformly scaling property is an ergodic-theoretical notion of self-similarity. Hochman proved the striking fact that generated distributions are always fractal distributions. The following result of K\"aenm\"aki, Sahlsten, and Shmerkin \cite[Theorem C]{KaenmakiSahlstenShmerkin2014a} is a converse to this.

\begin{theorem} \label{thm:usm}
  If $P$ is a fractal distribution, then there exists a uniformly scaling measure $\mu$ generating $P$.
\end{theorem}

Recall that if $P$ is an ergodic fractal distribution, then, by \eqref{eq:erg_usm}, $P$-almost every measure is uniformly scaling. Thus, by Theorems \ref{thm:compact} and \ref{thm:poulsen}, it suffices to show that the collection of fractal distributions satisfying the claim is closed. Let $(P_i)_i$ be a sequence of ergodic fractal distributions converging to $P$ and let $\mu_i$ be a uniformly scaling measure generating $P_i$. The proof is again based on the interplay between fractal distributions and CP processes. The rough idea to obtain a uniformly scaling measure generating $P$ is to splice the measures $\mu_i$ together. For the full proof, the reader is referred to \cite[\S 5]{KaenmakiSahlstenShmerkin2014a}.

\section{Geometry of measures}

Let $G(d,d-k)$ denote the set of all $(d-k)$-dimensional linear subspaces of $\R^d$. For $x \in \R^d$, $r>0$, $V \in G(d,d-k)$, and $0<\alpha\le 1$ define
\begin{equation*}
  X(x,r,V,\alpha) = \{ y \in B(x,r) : \dist(y-x,V) < \alpha|y-x| \}.
\end{equation*}
Conical density results aim to give conditions on a measure which guarantee that the cones $X(x,r,V,\alpha)$ contain a large portion of the mass from the surrounding ball $B(x,r)$ for certain proportion of scales. For example, a lower bound on some dimension often is such a condition. Recall that the \emph{lower local dimension} of a Radon measure $\mu$ at $x \in \R^d$ is
\begin{equation} \label{eq:lower_local_dim}
  \ldimloc(\mu,x) = \liminf_{r \downarrow 0} \frac{\log\mu(B(x,r))}{\log r}
\end{equation}
and the \emph{lower Hausdorff dimension} of $\mu$ is
\begin{align*}
  \ldimh(\mu) &= \essinf_{x \sim \mu} \ldimloc(\mu,x) \\
  &= \inf\{ \dimh(A) : A \subset \R^d \text{ is a Borel set with } \mu(A)>0 \}.
\end{align*}
Here $\dimh(A)$ is the Hausdorff dimension of the set $A \subset \R^d$. A measure $\mu$ is \emph{exact-dimensional} if the limit in \eqref{eq:lower_local_dim} exists and is $\mu$-almost everywhere constant. In this case, the common value is simply denoted by $\dim(\mu)$.

Intuitively, the local dimension of a measure should not be affected by the geometry of the measure on a density zero set of scales. Thus one could expect that tangent distributions should encode all information on dimensions.

\begin{definition}[Dimension of fractal distributions]
  The \emph{dimension of a fractal distribution} $P$ is
  \begin{equation*}
    \dim(P) = \int \dim(\mu) \dd P(\mu).
  \end{equation*}
\end{definition}

The dimension above is well defined by the fact that if $P$ is a fractal distribution, then $P$-almost every measure is exact-dimensional; see \cite[Lemma 1.18]{Hochman2010}. The dimension of fractal distributions has also other convenient properties. While the Hausdorff dimension is highly discontinuous on measures, the function $P \mapsto \dim(P)$ defined on the family of fractal distributions is continuous; see \cite[Lemma 3.20]{KaenmakiSahlstenShmerkin2014b}. The usefulness of the definition is manifested in the following result of Hochman \cite[Proposition 1.19]{Hochman2010}. Recall Theorem \ref{thm:hochman}.

\begin{theorem} \label{thm:fd_dim}
  If $\mu$ is a Radon measure, then
  \begin{equation*}
    \ldimloc(\mu,x) = \inf\{ \dim(P) : P \in \TD(\mu,x) \}
  \end{equation*}
  for $\mu$-almost all $x$. Furthermore, if $\mu$ is a uniformly scaling measure generating a fractal distribution $P$, then $\mu$ is exact-dimensional and $\dim(\mu) = \dim(P)$.
\end{theorem}

It turns out that tangent distributions are well suited to address problems concerning conical densities. The cones in question do not change under magnification and this allows to pass information between the original measure and its tangent distributions. Let
\begin{equation*}
  \AA_\eps := \{ \nu \in \MM_1 : \nu(X(0,1,V,\alpha)) \le \eps \text{ for some } V \in G(d,d-k) \}
\end{equation*}
for all $\eps \ge 0$. It is straightforward to see that $\AA_\eps$ is closed for all $\eps \ge 0$; see \cite[Lemma 4.2]{KaenmakiSahlstenShmerkin2014b}. The key observation is that
\begin{equation*}
  \AA_0 = \{ \nu \in \MM_1 : \spt(\nu) \cap X(0,1,V,\alpha) = \emptyset \text{ for some } V \in G(d,d-k) \},
\end{equation*}
where the defining property concerns only sets, is $S_t$ invariant.

The following conical density result is proved by K\"aenm\"aki, Sahlsten, and Shmerkin \cite[Proposition 4.3]{KaenmakiSahlstenShmerkin2014b}. Roughly speaking, it claims that if the dimension of the measure is large, then there are many scales in which the cones contain a relatively large portion of the mass. A slightly more precise version is that there exists $\eps > 0$ such that if $\ldimh(\mu) > k$, then for many scales $e^{-t}>0$ we have
\begin{equation*}
  \inf_{V \in G(d,d-k)} \frac{\mu(X(x,e^{-t},V,\alpha))}{\mu(B(x,e^{-t}))} > \eps
\end{equation*}
for $\mu$-almost all $x$. The precise formulation of the theorem is as follows.

\begin{theorem} \label{thm:conical}
  If $k \in \{ 1,\ldots,d-1 \}$, $k<s\le d$, and $0<\alpha\le 1$, then there exists $\eps > 0$ satisfying the following: For every Radon measure $\mu$ on $\R^d$ with $\ldimh(\mu) \ge s$ it holds that
  \begin{equation*} 
    \liminf_{T \to \infty} \la \mu \ra_{x,T}(\MM_1 \setminus \AA_\eps) \ge \frac{s-k}{d-k}
  \end{equation*}
  for $\mu$-almost all $x \in \R^d$.
\end{theorem}

The proof is based on showing that there cannot be ``too many'' rectifiable tangent measures. This means that, perhaps surprisingly, most of the known conical density results are, in some sense, a manifestation of rectifiability.

\begin{definition}[Rectifiability]
  A set $E \subset \R^d$ is called \emph{$k$-rectifiable} if there are countably many Lipschitz maps $f_i \colon \R^k \to \R^d$ so that
  \begin{equation*}
    \HH^k\Bigl( E \setminus \bigcup_i f_i(\R^k) \Bigr) = 0.
  \end{equation*}
\end{definition}

Here $\HH^k$ is the $k$-dimensional Hausdorff measure. Observe that a $k$-rectifiable set $E$ has $\dimh(E) \le k$. A sufficient condition for a set $E \subset \R^d$ to be $k$-rectifiable is that for every $x \in E$ there are $V \in G(d,d-k)$, $0 < \alpha < 1$, and $r > 0$ such that $E \cap X(x,r,V,\alpha) = \emptyset$; see \cite[Lemma 15.13]{Mattila1995}. Thus, if a fractal distribution $P$ satisfies $P(\AA_0) = 1$, then the quasi-Palm property implies that the support of $P$-almost every $\nu$ is $k$-rectifiable and hence $\dim(P) \le k$.

To prove Theorem \ref{thm:conical}, let $p,\delta>0$ be such that $p < (s-\delta-k)/(d-k) < (s-k)/(d-k)$. Suppose to the contrary that there is $0<\alpha\le 1$ so that for each $0<\eps<\eps(d,k,\alpha)$ there exists a Radon measure $\mu$ with $\ldimh(\mu) \ge s$ such that the claim fails to hold for $p$, that is,
\begin{equation*}
  \limsup_{T \to \infty}\,\la \mu \ra_{x,T}(\AA_\eps) > 1-p
\end{equation*}
on a set $E_\eps$ of positive $\mu$ measure. By Theorems \ref{thm:hochman} and \ref{thm:fd_dim}, we may assume that at points $x \in E_{\eps}$, all tangent distributions of $\mu$ are fractal distributions and
\begin{equation*}
  \inf\{\dim(P) : P \in \TD(\mu,x)\} = \ldimloc(\mu,x) > s-\delta.
\end{equation*}
Fix $x \in E_{\eps}$. For each $0<\eps<\eps(d,k,\alpha)$, as $\AA_\eps$ is closed, we find a tangent distribution $P_{\eps} \in \TD(\mu,x)$ so that $P_{\eps}(\AA_\eps) \ge 1-p$. Since the sets $\AA_\eps$ are also nested, we get
\begin{equation*}
  P(\AA_0) = \lim_{\eps \downarrow 0} P(\AA_\eps) \ge 1-p,
\end{equation*}
where $P$ is a weak limit of a sequence formed from $P_\eps$ as $\eps \downarrow 0$. Furthermore, since the collection of all fractal distributions is closed by Theorem \ref{thm:compact} and the dimension is continuous, the limit distribution $P$ is a fractal distribution with
\begin{equation*}
  \dim(P) \ge s-\delta.
\end{equation*}
Let $P_\omega$, $\omega \sim P$, be the ergodic components of $P$. By the invariance of $\AA_0$, we have $P_\omega(\AA_0) \in \{ 0,1 \}$ for $P$-almost all $\omega$. If $P_\omega(\AA_0) = 0$, then we use the trivial estimate $\dim(P_\omega) \le d$, and if $P_\omega(\AA_0) = 1$, then the rectifiability argument gives $\dim(P_\omega) \le k$. Since $P(\{ \omega : P_\omega(\AA_0) = 1 \}) = P(\AA_0) \ge 1-p$ we estimate
\begin{align*}
  s-\delta \le \dim(P) = \int \dim(P_\omega) \dd P(\omega)
  \le P(\AA_0)k + (1-P(\AA_0))d \le (1-p)k + pd
\end{align*}
yielding $p \ge (s-\delta-k)/(d-k)$. But this contradicts the choice of $\delta$. Thus the claim holds.

Relying on the existence of uniform scaling measures, we are able to study the sharpness of Theorem \ref{thm:conical}. The following result is proved by K\"aenm\"aki, Sahlsten, and Shmerkin \cite[Proposition 4.4]{KaenmakiSahlstenShmerkin2014b}.

\begin{theorem}
  If $k \in \{ 1,\ldots,d-1 \}$, $k < s \le d$, and $0<\alpha\le 1$, then there exists a Radon measure $\mu$ on $\R^d$ with $\dim(\mu) = s$ such that
  \begin{equation*}
    \lim_{T \to \infty} \la \mu \ra_{x,T}(\MM_1 \setminus \AA_\eps) = 
    \begin{cases}
      (s-k)/(d-k), &\text{if } 0 < \eps < \eps(d,k,\alpha), \\
      0,           &\text{if } \eps > \eps(d,k,\alpha),
    \end{cases}
  \end{equation*}
  for $\mu$-almost all $x \in \R^d$.
\end{theorem}

Here, for $k \in \{ 1,\ldots,d-1 \}$, $0 < \alpha \le 1$, and $V \in G(d,d-k)$, we have defined
\begin{equation*}
  \eps(d,k,\alpha) := \frac{\LL^d(X(0,1,V,\alpha))}{\LL^d(B(0,1))}.
\end{equation*}
It follows from the rotational invariance of the Lebesgue measure $\LL^d$ that $\eps(d,k,\alpha)$ does not depend on the choice of $V$.

The measure $\mu$ above is just a uniform scaling measure generating
\begin{equation*}
  P = \frac{s-k}{d-k} \delta_\LL + \Bigl( 1-\frac{s-k}{d-k} \Bigr) \delta_\HH,
\end{equation*}
where $\LL$ is the normalization of $\LL^d|_{B_1}$ and $\HH$ is the normalization of $\HH^k|_{W \cap B_1}$ for a fixed $W \in G(d,k)$. Since $P$ is a convex combination of two fractal distributions, it is a fractal distribution. The existence of $\mu$ is guaranteed by Theorem \ref{thm:usm}. Recalling Theorem \ref{thm:fd_dim}, we see that $\mu$ is exact-dimensional and
\begin{equation*}
  \dim(\mu) = \dim(P) = \frac{s-k}{d-k}\,d + \Bigl( 1 - \frac{s-k}{d-k} \Bigr)k = s.
\end{equation*}
The goal is to verify that $\mu$ has the claimed properties.

Fix $0 < \eps < \eps(d,k,\alpha)$. Since $\LL(X(0,1,V,\alpha)) = \eps(d,k,\alpha) > \eps$ for all $V \in G(d,d-k)$ and $\HH(X(0,1,W^\bot,\alpha)) = 0$ we have $P(\MM_1 \setminus \AA_\eps) = (s-k)/(d-k)$. Thus, by the weak convergence, it follows that
\begin{equation*}
  \lim_{T \to \infty} \la \mu \ra_{x,T}(\MM_1 \setminus \AA_\eps) = \frac{s-k}{d-k}.
\end{equation*}
In the case $\eps > \eps(d,k,\alpha)$ we can reason similarly.

\end{document}